\documentclass[12pt]{article}
\usepackage{amsmath,amssymb}
\usepackage{graphicx}

\newtheorem{theorem}{Theorem}[section]

\newtheorem{corollary}{Corollary}[section]

\newtheorem{Rem}{Remark}[section]

\numberwithin{equation}{section}

\def\bbe{\mathbb E}

\def\bbp{\mathbb P}
\def\bbr{\mathbb R}

\begin{document}

\title{Asymptotics for the length of the longest increasing subsequence of a binary Markov 
random word}
\author{Christian Houdr\'e \thanks{Georgia Institute of Technology, 
School of Mathematics, Atlanta, Georgia, 30332-0160, USA, 
houdre@math.gatech.edu. {Research supported in part by the NSA Grant H98230-09-1-0017.}} 
\and Trevis J. Litherland \thanks{Georgia Institute of Technology, 
School of Mathematics, Atlanta, Georgia, 30332-0160, USA, 
trevisl@math.gatech.edu} }

\maketitle

\vspace{0.5cm}

\begin{abstract}
\noindent Let $(X_n)_{n \ge 0}$ 
be an irreducible, aperiodic, and homogeneous
binary Markov chain and let $LI_n$ be the length of the 
longest (weakly) increasing subsequence of $(X_k)_{1\le k \le n}$.  
Using combinatorial constructions and 
weak invariance principles, we present elementary 
arguments leading to a new proof that (after proper centering and scaling) 
the limiting law of $LI_n$ is the maximal eigenvalue of a $2\times 2$ 
Gaussian random matrix.  In fact, the limiting shape of the RSK Young diagrams 
associated with the binary Markov random word is the spectrum of this 
random matrix.  
\end{abstract}





\section{Introduction}  

The identification of the limiting distribution of 
the length of the longest increasing
subsequence  of a random permutation or of a random word has attracted 
a lot of interest in the past decade, in particular in light of 
its connections with random matrices 
(see \cite{BDJ1, BDJ2, BDJ2Add}, \cite{Ba}, \cite{BOO}, \cite{GTW}, \cite{HX}, 
\cite{ITW1,ITW2}, \cite{Jo}, \cite{Ok}, \cite{TW}).  For random words, 
both the iid uniform and non-uniform settings are understood, leading  
respectively to the maximal eigenvalue of a traceless (or generalized traceless) 
element of the Gaussian Unitary Ensemble (GUE) as limiting laws of $LI_n$.  
In a dependent framework, Kuperberg \cite{Ku} conjectured that
if the word is generated by an irreducible, doubly-stochastic, cyclic, Markov chain with state  
space an ordered $m$-letter alphabet, 
then the limiting distribution of the length $LI_n$ is still that of the 
maximal eigenvalue of a traceless $m \times m$ element of the GUE.
More generally, the conjecture asserts that the shape of the Robinson-Schensted-Knuth (RSK) Young diagrams 
associated with the Markovian random word is that of 
the joint distribution of the eigenvalues
of a traceless $m \times m$ element of the GUE.  
For $m=2$, Chistyakov and G\"otze \cite{ChG} positively answered this 
conjecture, and in the present 
paper this result is 
rederived in an elementary way.

The precise class of homogeneous Markov chains 
with which Kuperberg's conjecture is concerned is more specific than
the ones we shall study.  The irreducibility of
the chain is a basic property we certainly must demand:
each letter has to occur at some point following the
occurrence of any given letter.  The cyclic (also called {\it circulant})
criterion, {\it i.e.}, the Markov transition
matrix $P$ has entries satisfying $p_{i,j} = p_{i+1,j+1}$,
for $1 \le i,j \le m$ (where $m+1 = 1$), ensures 
a uniform stationary distribution.  

Let us also note that Kuperberg implicitly assumes 
the Markov chain to also be aperiodic.  Indeed, the simple
$2$-state Markov chain for the letters $\alpha_1$ and $\alpha_2$
described by $\bbp(X_{n+1} = \alpha_i| X_n = \alpha_j) = 1$ for $i \ne j$,
produces a sequence of alternating letters, so that
$LI_n$ is always either $n/2$ or $n/2 + 1$, for $n$ even,
and $(n+1)/2$, for $n$ odd, and so has a degenerate
limiting distribution.  Even though this Markov chain
is irreducible and cyclic,
it is periodic.

By the end of this introduction, the reader might certainly 
have wondered how the binary 
results do get modified for ordered alphabets of arbitrary fixed size $m$.  
As shown in \cite{HL2}, for $m=3$, Kuperberg's conjecture is indeed true. 
However, for $m\ge4$, this is no longer the case; and some, 
but not all, cyclic Markov chains 
lead to a limiting law as in the iid uniform case.

\section{Combinatorics}

As in \cite{HL}, one can express 
$LI_n$ in a combinatorial manner.
For convenience, this short section recapitulates that development.\\

Let $(X_n)_{n\ge1}$ consist of a sequence of values taken from an
$m$-letter ordered alphabet, 
${\cal A}_m = \{\alpha_1 < \alpha_2 < \cdots < \alpha_m\}$. 
Let $a^r_k$ be the number of occurrences
of $\alpha_r$ among $(X_i)_{1\le i \le k}$.  
Each increasing subsequence of $(X_i)_{1\le i \le k}$ consists simply of 
consecutive identical values, 
with these values forming an increasing subsequence of $\alpha_r$.
Moreover, the number of occurrences of $\alpha_r\in \{\alpha_1,
\dots,\alpha_m\}$ among
$(X_i)_{k+1 \le i \le \ell}$,
where $1 \le k < \ell \le n$, is simply
$a^r_{\ell}-a^r_k$. The length of the longest increasing subsequence 
of $X_1,X_2,\dots, X_n$ is thus given by

\begin{equation}\label{item1}
LI_n=\max_{\stackrel{\scriptstyle 0\le k_1\le\cdots}{\le k_{m-1}\le n}}
     [(a^1_{k_1}-a^1_0)+(a^2_{k_2}-a^2_{k_1})+\cdots +
     (a^m_n-a^m_{k_{m-1}})],
\end{equation}

\noindent {\it i.e.},

\begin{equation}\label{item2}
LI_n=\max_{\stackrel{\scriptstyle 0\le k_1\le\cdots}{\le k_{m-1}\le n}}
     [(a^1_{k_1}-a^2_{k_1})+(a^2_{k_2}-a^3_{k_2})+\cdots +
     (a^{m-1}_{k_{m-1}}-a^m_{k_{m-1}})+a^m_n],
\end{equation}

\noindent where $a^r_0=0$.

For $i = 1, \dots ,n$ and $r = 1, \dots ,{m-1}$, let 

\begin{equation}\label{item3}
Z^r_i=	\begin{cases} 	1, &\text{if $X_i=\alpha_r,$}\\
			-1, & \text{if $X_i=\alpha_{r+1},$}\\
			0, & \text{otherwise,}
	\end{cases}
\end{equation}

\noindent and let $S^r_k=\sum^k_{i=1}Z^r_i$, $k = 1, \dots ,n$, with also $S^r_0=0$. 
Then clearly $S^r_k = a^r_k - a^{r+1}_k$.  Hence,

\begin{equation}\label{item4}
LI_n=\max_{\stackrel{\scriptstyle 0\le k_1\le\cdots}{\le k_{m-1}\le n}}
     \{S^1_{k_1}+S^2_{k_2}+ \cdots + S^{m-1}_{k_{m-1}}+a^m_n\}.
\end{equation}

By the telescoping nature of the sum $\sum_{k=r}^{m-1}S_n^k = 
\sum_{k=r}^{m-1}(a_n^k - a_n^{k+1})$,
we find that, for each $1 \le r \le m-1$,
$a_n^r = a_n^m + \sum_{k=r}^{m-1}S_n^k$.
Since $a^1_k, \dots ,a^m_k$ must evidently sum up 
to $k$, we have

\begin{align*}
n &= \sum^m_{r=1}a^r_n \\
&= \sum^{m-1}_{r=1} \left(a_n^m + \sum_{k=r}^{m-1}S_n^k\right) + a^m_n\\
&= \sum^{m-1}_{r=1}rS^r_n + ma^m_n.\\
\end{align*}

Solving for $a^m_n$ gives us
$$a^m_n=\frac nm- \frac1m \sum^{m-1}_{r=1} rS^r_n.$$

Substituting into \eqref{item4}, we finally obtain

\begin{equation}\label{item5}
 LI_n=\frac nm-\frac1m\sum^{m-1}_{r=1} rS^r_n+
      \max_{\stackrel{\scriptstyle 0\le k_1\le\cdots}{\le k_{m-1}\le n}}
      \{S^1_{k_1}+S^2_{k_2}+ \cdots + S^{m-1}_{k_{m-1}}\}.
\end{equation}

As emphasized in \cite{HL}, \eqref{item5} is of a {\it purely
combinatorial nature or, in more probabilistic terms, is of a pathwise nature}.  
We now proceed to analyze \eqref{item5} for a binary Markovian sequence.

\section{Binary Markovian Alphabet}

In the context of binary Markovian alphabets,  
$(X_n)_{n\ge0}$ is described by the following
transition probabilities between the two states
(which we identify with the two letters 
$\alpha_1$ and $\alpha_2$):
$\bbp(X_{n+1} = \alpha_2 | X_n = \alpha_1) = a$ and
$\bbp(X_{n+1} = \alpha_1 | X_n = \alpha_2) = b$,
where $0 < a + b < 2$.  We later examine the degenerate
cases $a = b = 0$ and $a = b = 1$.
In keeping with the common usage within the Markov chain
literature, we begin our sequence at $n=0$, although
our focus will be on $n \ge 1$.
Denoting by $(p_n^1,p_n^2)$
the vector describing the probability distribution
on $\{\alpha_1,\alpha_2\}$ at time $n$, 
we have

\begin{equation}\label{item5a}
  \begin{pmatrix} p_{n+1}^1, p_{n+1}^2 \end{pmatrix}
= \begin{pmatrix} p_{n}^1, p_{n}^2   \end{pmatrix}
  \begin{pmatrix} 1-a & a\\ b & 1-b   \end{pmatrix}.
\end{equation}

The eigenvalues of the
matrix in \eqref{item5a} are $\lambda_1 = 1$ and 
$-1 < \lambda_2 = 1 - a - b < 1$, 
with respective left eigenvectors
$(\pi_1, \pi_2) = (b/(a+b),a/(a+b))$ and
$(1,-1)$.  Moreover, $(\pi_1, \pi_2)$ is also the
stationary distribution.  Given any initial
distribution $(p_0^1,p_0^2)$, we find that

\begin{equation}\label{item5b}
\begin{pmatrix} p_{n}^1, p_{n}^2 \end{pmatrix}
=  \begin{pmatrix} \pi_1, \pi_2  \end{pmatrix}
+ \lambda_2^{n}\frac{a p_0^1-b p_0^2}{a+b} 
  \begin{pmatrix} 1,-1\end{pmatrix}
  \rightarrow \begin{pmatrix} \pi_1, \pi_2  \end{pmatrix},
\end{equation}

\noindent as $n \rightarrow \infty$, since $|\lambda_2| < 1$.\\

Our goal is now to use these probabilistic
expressions to describe the random variables
$Z_k^1$ and $S_k^1$ defined in the previous section.  
(We retain the redundant
superscript ``$1$'' in $Z_k^1$ and $S_k^1$ in the 
interest of uniformity.)

Setting $\beta = a p_0^1-b p_0^2$,
we easily find that 

\begin{align}\label{item5c}
\bbe Z_k^1 &= (+1)\left(\pi_1 + \frac{\beta}{a+b}\lambda_2^{k} \right)
            + (-1)\left(\pi_2 - \frac{\beta}{a+b}\lambda_2^{k} \right)\nonumber\\
	   &= \frac{b-a}{a+b}  + 2 \frac{\beta}{a+b}\lambda_2^{k},
\end{align}

\noindent for each $1 \le k \le n$.  Thus,

\begin{equation}\label{item5d}
\bbe S_k^1 = \frac{b-a}{a+b}k  
           + 2 \left(\frac{\beta \lambda_2}{a+b}\right) \left( \frac{1-\lambda_2^k}{1-\lambda_2} \right),
\end{equation}

\noindent and so $\bbe S_k^1/k \rightarrow (b-a)/(a+b)$,
as $k \rightarrow \infty$.

Turning to the second moments of 
$Z_k^1$ and $S_k^1$, first note
that $\bbe (Z_k^1)^2 = 1$, since
$(Z_k^1)^2 = 1$ a.s.
Next, we consider $\bbe Z_k^1 Z_{\ell}^1$,
for $k < \ell$.  Using
the Markovian structure of $(X_n)_{n\ge0}$,
it quickly follows that

\begin{align}\label{item5e}
&\bbp((X_k,X_{\ell}) = (x_k,x_{\ell}))\nonumber\\
 &\qquad =
\begin{cases}
  \left(\pi_1 + \lambda_2^{\ell-k}\frac{a}{a+b}\right)\left(\pi_1 + \lambda_2^{k}\frac{\beta}{a+b}\right),
     &\text{if $(x_k,x_{\ell})=(\alpha_1,\alpha_1)$},\\
  \left(\pi_1 - \lambda_2^{\ell-k}\frac{b}{a+b}\right)\left(\pi_2 - \lambda_2^{k}\frac{\beta}{a+b}\right),
     &\text{if $(x_k,x_{\ell})=(\alpha_1,\alpha_2)$},\\
  \left(\pi_2 - \lambda_2^{\ell-k}\frac{a}{a+b}\right)\left(\pi_1 + \lambda_2^{k}\frac{\beta}{a+b}\right),
     &\text{if $(x_k,x_{\ell})=(\alpha_2,\alpha_1)$},\\
  \left(\pi_2 + \lambda_2^{\ell-k}\frac{b}{a+b}\right)\left(\pi_2 - \lambda_2^{k}\frac{\beta}{a+b}\right),
     &\text{if $(x_k,x_{\ell})=(\alpha_2,\alpha_2)$}.
\end{cases}
\end{align}

For simplicity, we will henceforth 
assume that our initial 
distribution is the stationary one, {\it i.e.},
$(p_0^1,p_0^2) = (\pi_1, \pi_2)$.
(this assumption can be dropped as explained 
in the Concluding Remarks of \cite{HL2}).  
Under this assumption, $\beta = 0$, 
$\bbe S_k^1 = k\mu$, where 
$\mu = \bbe Z_k^1 = (b-a)/(a+b)$, and 
\eqref{item5e} simplifies to

\begin{align}\label{item5f}
&\bbp((X_k,X_{\ell}) = (x_k,x_{\ell}))\nonumber\\
 &\qquad =
\begin{cases}
  \left(\pi_1 + \lambda_2^{\ell-k}\frac{a}{a+b}\right)\pi_1,
     &\text{if $(x_k,x_{\ell})=(\alpha_1,\alpha_1)$},\\
  \left(\pi_1 - \lambda_2^{\ell-k}\frac{b}{a+b}\right)\pi_2,
     &\text{if $(x_k,x_{\ell})=(\alpha_1,\alpha_2)$},\\
  \left(\pi_2 - \lambda_2^{\ell-k}\frac{a}{a+b}\right)\pi_1,
     &\text{if $(x_k,x_{\ell})=(\alpha_2,\alpha_1)$},\\
  \left(\pi_2 + \lambda_2^{\ell-k}\frac{b}{a+b}\right)\pi_2,
     &\text{if $(x_k,x_{\ell})=(\alpha_2,\alpha_2)$}.
\end{cases}
\end{align}

We can now compute $\bbe Z_k^1 Z_{\ell}^1$:

{\allowdisplaybreaks
\begin{align}\label{item5ga}
\bbe Z_k^1 Z_{\ell}^1 
&= \bbp(Z_k^1 Z_{\ell}^1 = +1) - \bbp(Z_k^1 Z_{\ell}^1 = -1)\nonumber\\
&= \bbp( (X_k,X_{\ell}) \in \{(\alpha_1,\alpha_1),(\alpha_2,\alpha_2) \})\nonumber\\
  &\qquad\qquad - \bbp( (X_k,X_{\ell}) \in \{(\alpha_1,\alpha_2),(\alpha_2,\alpha_1) \})\nonumber\\
&= \left( \pi_1^2 + \lambda_2^{\ell-k}\frac{a}{a+b}\pi_1 + \pi_2^2 + \lambda_2^{\ell-k}\frac{b}{a+b}\pi_2\right)\nonumber\\
  &\qquad\qquad -\left( \pi_1 \pi_2 - \lambda_2^{\ell-k}\frac{b}{a+b}\pi_2 + \pi_1 \pi_2  - \lambda_2^{\ell-k}\frac{a}{a+b}\pi_1\right)\nonumber\\
&= \left(\pi_1^2 + \pi_2^2 + \frac{2ab}{(a+b)^2}\lambda_2^{\ell-k}\right)
  -\left(2\pi_1 \pi_2 - \frac{2ab}{(a+b)^2}\lambda_2^{\ell-k}\right)\nonumber\\
&= \frac{(b-a)^2}{(a+b)^2} + \frac{4ab}{(a+b)^2}\lambda_2^{\ell-k}.
\end{align}
}

Hence, recalling that $\beta = 0$,

\begin{align}\label{item5gb}
\sigma^2 := \mbox{Var} Z_k^1 
&= 1 - \left(\frac{b-a}{a+b}\right)^2\nonumber\\
&= \frac{4ab}{(a+b)^2},
\end{align}

\noindent for all $k\ge 1$,
and, for $k < \ell$, the covariance
of $Z_k^1$ and $Z_{\ell}^1$ is

\begin{align}\label{item5h}
\mbox{Cov} (Z_k^1,Z_{\ell}^1) &= \frac{(b-a)^2}{(a+b)^2} + \sigma^2\lambda_2^{\ell-k} - \left(\frac{b-a}{a+b}\right)^2 = \sigma^2\lambda_2^{\ell-k}.
\end{align}

Proceeding to the covariance structure
of $S_k^1$, we first find that

\begin{align}\label{item5i}
\mbox{Var}S_k^1 &= \sum_{j=1}^k \mbox{Var} Z_j^1 + 2 \sum_{j<\ell}\mbox{Cov}(Z_j^1,Z_l^1)\nonumber\\
&= \sigma^2k + 2\sigma^2\sum_{j<\ell}\lambda_2^{\ell-j}\nonumber\\
&= \sigma^2k 
  + 2\sigma^2\left(\frac{\lambda_2^{k+1} - k\lambda_2^2 + (k-1)\lambda_2}{(1-\lambda_2)^2}\right)\nonumber\\
&= \sigma^2\left(\frac{1+\lambda_2}{1-\lambda_2}\right)k
 + 2\sigma^2 \left(\frac{\lambda_2(\lambda_2^k-1)}{(1-\lambda_2)^2}\right).
\end{align}

Next, for $k < \ell$, and using \eqref{item5h} and \eqref{item5i},
the covariance of $S_k^1$ and $S_{\ell}^1$ is
given by

{\allowdisplaybreaks
\begin{align}
\mbox{Cov} (S_k^1,S_{\ell}^1) 
  &= \sum_{i=1}^k \sum_{j=1}^{\ell} \mbox{Cov} (Z_i^1,Z_j^1)\nonumber\\
  &= \sum_{i=1}^k \mbox{Var} Z_i^1 + 2\sum_{i<j<k} \mbox{Cov} (Z_i^1,Z_j^1) 
      + \sum_{i=1}^k \sum_{j=k+1}^{\ell} \mbox{Cov} (Z_i^1,Z_j^1)\nonumber\\ 
  &= \mbox{Var} S_k^1 + \sum_{i=1}^k \sum_{j=k+1}^{\ell} \mbox{Cov} (Z_i^1,Z_j^1)\nonumber\\ 
  &= \mbox{Var} S_k^1 
    +  \sigma^2\left(\frac{\lambda_2(1-\lambda_2^k)(1-\lambda_2^{\ell-k})}{(1-\lambda_2)^2}\right)\nonumber\\ 
  &= \sigma^2\left( \left(\frac{1+\lambda_2}{1-\lambda_2}\right)k
    - \frac{\lambda_2(1-\lambda_2^k)(1+\lambda_2^{\ell-k})}{(1-\lambda_2)^2}\right).\label{item5jb}
\end{align}
}

From \eqref{item5i} and \eqref{item5jb} we see that,
as $k \rightarrow \infty$,

\begin{equation}\label{item5k}
\frac{\mbox{Var}S_k^1}{k} \rightarrow \sigma^2\left(\frac{1+\lambda_2}{1-\lambda_2}\right),
\end{equation}

\noindent and, moreover,
as $k \wedge \ell \rightarrow \infty$,

\begin{equation}\label{item5l}
\frac{\mbox{Cov} (S_k^1,S_{\ell}^1)}{(k \wedge \ell)} \rightarrow \sigma^2\left(\frac{1+\lambda_2}{1-\lambda_2}\right).
\end{equation}

\noindent When $a = b$, $\bbe S_k^1 = 0$, and in \eqref{item5k}
the asymptotic variance becomes

\begin{align*}
 \frac{\mbox{Var}S_k^1}{k} 
 &\rightarrow  \frac{4a^2}{(2a)^2}\left(\frac{1+(1-2a)}{1-(1-2a)}\right)\nonumber\\
 &= \frac{1}{a} - 1.
\end{align*}

For $a$ small, we have a ''lazy'' Markov chain,
that is, a Markov chain which tends to remain
in a given state for long periods of time.
In this regime, the random variable $S_k^1$ has long
periods of increase followed by long periods
of decrease.  In this way, linear asymptotics
of the variance with large constants
occur.  If, on the other hand, $a$ is close to $1$,
the Markov chain rapidly shifts back and forth
between $\alpha_1$ and $\alpha_2$, and so
the constant associated with the linearly 
increasing variance of  $S_k^1$ is small.

As in \cite{HL}, Brownian functionals play a central r\^ole in
describing the limiting distribution of $LI_n$.  
To move towards a Brownian functional
expression for the limiting law of $LI_n$, 
define the polygonal function

\begin{equation}\label{item5n}
\hat B_n(t)=\frac{S^1_{[nt]} - [nt]\mu}{\sigma\sqrt{n(1+\lambda_2)/(1-\lambda_2)}} 
+\frac{(nt-[nt])(Z^1_{[nt]+1} - \mu)}{\sigma\sqrt{n(1+\lambda_2)/(1-\lambda_2)}}, 
\end{equation}

\noindent for $0 \le t \le 1$.  In our finite-state,
irreducible, aperiodic, stationary Markov chain setting, 
we may conclude that $\hat B_n \Rightarrow B$, as desired.
(See, for example, the more general settings for
Gordin's martingale approach to dependent
invariance principles, and the stationary ergodic
invariance principle found in Theorem 19.1 of \cite{Bill}.)

Turning now to $LI_n$, we see that
for the present $2$-letter situation,
\eqref{item5} simply becomes

$$LI_n = \frac{n}{2} - \frac{1}{2}S^1_n + \max_{1 \le k \le n} S^1_k.$$

To find the limiting distribution of $LI_n$ from this
expression, recall that $\pi_1 = b/(a+b)$,
$\pi_2 = a/(a+b)$, $\mu = \pi_1 - \pi_2 = (b-a)/(a+b)$,
$\sigma^2 = 4ab/(a+b)^2$, 
and that $\lambda_2 = 1-a-b$.
Define $\pi_{max} = \max\{\pi_1,\pi_2\}$ and
$\tilde{\sigma}^2 = \sigma^2(1+\lambda_2)/(1-\lambda_2)$.
Rewriting \eqref{item5n} as


$$ \hat B_n(t)=\frac{S^1_{[nt]} - [nt]\mu}{\tilde{\sigma}\sqrt{n}} 
+\frac{(nt-[nt])(Z^1_{[nt]+1} - \mu)}{\tilde{\sigma}\sqrt{n}}, $$
\
\noindent $LI_n$ becomes

\begin{align}\label{item5q}
LI_n &= \frac{n}{2}  -\frac12 \left(\tilde{\sigma}\sqrt{n} \hat B_n(1) + \mu n\right) 
+ \max_{0\le t \le 1} \left( \tilde{\sigma}\sqrt{n} \hat B_n(t) + \mu nt\right) \nonumber\\
&= n\pi_2  -\frac12 \left(\tilde{\sigma}\sqrt{n} \hat B_n(1)\right) + \max_{0\le t \le 1} 
\left( \tilde{\sigma}\sqrt{n} \hat B_n(t) + (\pi_1-\pi_2) nt\right) \nonumber\\
&= n\pi_{max}  -\frac12 \left(\tilde{\sigma}\sqrt{n} \hat B_n(1)\right) \nonumber\\
&\qquad + \max_{0\le t \le 1} \left( \tilde{\sigma}\sqrt{n} \hat B_n(t) + (\pi_1-\pi_2) nt - (\pi_{max} - \pi_2)n\right).
\end{align}

\noindent This immediately gives

\begin{align}\label{item5qa}
\frac{LI_n  - n\pi_{max}}{\tilde{\sigma}\sqrt{n}} &= -\frac12 \hat B_n(1) \nonumber\\
&+ \max_{0\le t \le 1} \left(  \hat B_n(t) + \frac{\sqrt{n}}{\tilde{\sigma}}((\pi_1-\pi_2)t - (\pi_{max} - \pi_2)) \right).
\end{align}

Let us examine \eqref{item5qa} on a case-by-case basis.
First, if $\pi_{max} = \pi_1 = \pi_2 = 1/2$, {\it i.e.},
if $a = b$, then
$\sigma = 1$ and $\tilde{\sigma} = (1-a)/a$, and
so \eqref{item5qa} becomes 

\begin{align}\label{item5qb}
\frac{LI_n  - n/2}{\sqrt{(1-a)n/a}} &= -\frac12 \hat B_n(1) + \max_{0\le t \le 1}  \hat B_n(t).
\end{align}

\noindent Then, by the Invariance Principle and the Continuous Mapping Theorem,

\begin{equation}\label{item5qba}
\frac{LI_n  - n/2}{\sqrt{(1-a)n/a}} \Rightarrow -\frac12 B(1) + \max_{0\le t \le 1}  B(t).
\end{equation}

Next, if $\pi_{max} = \pi_2 > \pi_1$, \eqref{item5qa} becomes

\begin{align}\label{item5qd}
\frac{LI_n  - n\pi_{max}}{\tilde{\sigma}\sqrt{n}} &= -\frac12 \hat B_n(1) \nonumber\\
&\qquad + \max_{0\le t \le 1} \left(  \hat B_n(t) - \frac{\sqrt{n}}{\tilde{\sigma}}(\pi_{max} - \pi_1)t \right).
\end{align}

On the other hand, if $\pi_{max} = \pi_1 > \pi_2$, \eqref{item5qa} becomes

\begin{align}\label{item5qc}
\frac{LI_n  - n\pi_{max}}{\tilde{\sigma}\sqrt{n}} 
&= -\frac12 \hat B_n(1) \nonumber\\
&\qquad + \max_{0\le t \le 1} \left(  \hat B_n(t) - 
\frac{\sqrt{n}}{\tilde{\sigma}}(\pi_{max} - \pi_2)(1-t) \right)\nonumber\\
&= \frac12 \hat B_n(1) \nonumber\\
&\qquad + \max_{0\le t \le 1} \left(  \hat B_n(t) - 
\hat B_n(1)- \frac{\sqrt{n}}{\tilde{\sigma}}(\pi_{max} - \pi_2)(1-t) \right).
\end{align}

In both \eqref{item5qd} and \eqref{item5qc} we have a term in
our maximal functional which is linear in $t$ or $1-t$,
with a negative slope.  We now show,
in an elementary fashion, that in both cases,
as $n \rightarrow \infty$, the maximal functional goes
to zero in probability.  

Consider first \eqref{item5qd}.  Let
$c_n = \sqrt{n}(\pi_{max} - \pi_1)/\tilde{\sigma} > 0,$
and for any $c > 0$, let 
$M_c = \max_{0\le t \le 1} (B(t) - ct)$,
where $B(t)$ is a standard Brownian motion.
Now for $n$ large enough,

$$\hat B_n(t) - c t \ge \hat B_n(t) - c_nt$$

\noindent a.s., for all $0 \le t \le 1$.
Then for any $z > 0$, and $n$ large enough, 

\begin{align}\label{item5qe}
\bbp( \max_{0 \le t \le 1} (\hat B_n(t) - c_n t) > z ) 
&\le \bbp( \max_{0 \le t \le 1} (\hat B_n(t) - c t) > z ),
\end{align}

\noindent and so by the Invariance Principle and the
Continuous Mapping Theorem,

\begin{align}\label{item5qf}
\limsup_{n \rightarrow \infty} \bbp( \max_{0 \le t \le 1} (\hat B_n(t) - c_n t) > z ) 
&\le \lim_{n \rightarrow \infty} \bbp( \max_{0 \le t \le 1} (\hat B_n(t) - c t) > z )\nonumber\\
&= \bbp( M_c > z ).
\end{align}

Now, as is well-known,
$\bbp( M_c > z ) \rightarrow 0$ as $c \rightarrow \infty$.
One can confirm this intuitive fact with the following simple argument.
For $z > 0$, $c > 0$, and $0 < \varepsilon < 1$, we have that

{\allowdisplaybreaks
\begin{align}\label{item5qfa}
\bbp ( M_c > z) &\le \bbp (\max_{0 \le t \le \varepsilon} (B(t) - c t) > z ) 
+ \bbp (\max_{\varepsilon < t \le 1} (B(t) - c t) > z )\nonumber\\
&\le  \bbp (\max_{0 \le t \le \varepsilon} B(t) > z ) + 
\bbp (\max_{\varepsilon < t \le 1} (B(t) - c \varepsilon) > z )\nonumber\\
&\le  \bbp (\max_{0 \le t \le \varepsilon} B(t) > z ) + \bbp (\max_{0 < t \le 1} B(t) >  c \varepsilon + z )\nonumber\\
&= 2\left(1 - \Phi\left(\frac{z}{\sqrt{\varepsilon}}\right) \right) + 2\left(1 - \Phi(c\varepsilon + z) \right).
\end{align}
}

\noindent But, as $c$ and $\varepsilon$ are arbitrary, we can
first take the limsup of \eqref{item5qfa} as $c \rightarrow \infty$,
and then let $\varepsilon \rightarrow 0$, 
proving the claim.

We have thus shown that 

$$\limsup_{n \rightarrow \infty} \bbp( \max_{0 \le t \le 1} (\hat B_n(t) - c_n t) > z ) \le 0,$$

\noindent and since the functional clearly is equal to zero when
$t = 0$, we have

\begin{equation}\label{item5qg}
\max_{0 \le t \le 1} (\hat B_n(t) - c_n t) \stackrel{\bbp}{\rightarrow} 0,
\end{equation}

\noindent as $n \rightarrow \infty$.
Thus, by the Continuous Mapping Theorem, and the Converging Together Lemma,
we obtain the weak convergence result

\begin{equation}\label{item5qh}
\frac{LI_n  - n\pi_{max}}{\tilde{\sigma}\sqrt{n}} \Rightarrow -\frac12 B(1).
\end{equation}

Lastly, consider \eqref{item5qc}.  
Here we need simply note the 
following equality in law, which
follows from the stationary and Markovian
nature of the underlying sequence $(X_n)_{n\ge0}$:

{\allowdisplaybreaks
\begin{align}\label{item5qi}
\hat B_n(t) - \hat B_n(1)&- \frac{\sqrt{n}}{\tilde{\sigma}}(\pi_{max} - \pi_2)(1-t)\nonumber\\
&\stackrel{\cal{L}}{=} -\hat B_n(1-t)- \frac{\sqrt{n}}{\tilde{\sigma}}(\pi_{max} - \pi_2)(1-t),
\end{align}
}

\noindent 
for $t = 0,1/n,\dots,(n-1)/n,1$.
With a change of variables
$(u = 1 - t)$, and noting that $B(t)$
and $-B(t)$ are equal in law,
our previous convergence result
\eqref{item5qg} implies that

\begin{align}\label{item5qj}
\max_{0 \le t \le 1} (\hat B_n(t) - \hat B_n(1)- c_n(1-t))
\stackrel{\cal{L}}{=}  \max_{0 \le u \le 1} (-\hat B_n(u) - c_n u)
\stackrel{\bbp}{\rightarrow} 0,
\end{align}

\noindent as $n \rightarrow \infty$.
Our limiting functional is thus of the form

\begin{equation}\label{item5ql}
\frac{LI_n  - n\pi_{max}}{\tilde{\sigma}\sqrt{n}} \Rightarrow \frac12 B(1).
\end{equation}

\noindent Since $B(1)$ is simply a standard normal random variable,
the different signs in \eqref{item5qh} and \eqref{item5ql}
are inconsequential.

Finally, consider the degenerate cases.
If either $a=0$ or $b = 0$, 
then the sequence $(X_n)_{n \ge 0}$ will be a.s.~constant,
regardless of the starting state, and so
$LI_n \sim n$.
On the other hand, if $a=b=1$, then the sequence
oscillates back and forth between $\alpha_1$
and $\alpha_2$, so that $LI_n \sim n/2$.
Combining these trivial cases with the previous development,
gives:

\begin{theorem}\label{thm1}
Let $(X_n)_{n \ge 0}$ be a $2$-state Markov chain,
with
$\bbp(X_{n+1} = \alpha_2|$\\
$X_n = \alpha_1) = a$ and
$\bbp(X_{n+1} = \alpha_1 | X_n = \alpha_2) = b$.
Let the law of $X_0$ be the invariant distribution
$(\pi_1,\pi_2) = (b/(a+b),a/(a+b))$, for $0 < a + b \le 2$, and be 
$(\pi_1,\pi_2) = (1,0)$, for  $a = b = 0$.
Then, for $a=b > 0$,

\begin{equation}\label{item5r}
\frac{LI_n-n/2}{\sqrt {n}} 
\Rightarrow  \sqrt{\frac{1-a}{a}}\left(-\frac12 B(1)+ \max_{0\le t \le 1}B(t)\right),
\end{equation}

\noindent 
where $(B(t))_{0\le t \le 1}$ is a standard Brownian motion.  
For $a \ne b$ or $a=b=0$, 

\begin{equation}\label{item5s}
\frac{LI_n  - n\pi_{max}}{\sqrt{n}} \Rightarrow N(0,\tilde{\sigma}^2/4),
\end{equation}

\noindent with $\pi_{max} = \max\{\pi_1,\pi_2\}$, and 
where $N(0,\tilde{\sigma}^2/4)$ 
is a centered normal
random variable with 
variance $\tilde{\sigma}^2/4 = ab(2-a-b)/(a+b)^3$, for $a \ne b$,
and $\tilde{\sigma}^2=0$, for $a=b=0$.
(If $a=b=1$, or $\tilde{\sigma}^2 = 0$, then the distributions in \eqref{item5r} and \eqref{item5s},
respectively, are understood to be degenerate at the origin.)
\end{theorem}

To extend this result to the entire RSK Young diagrams, let us introduce 
the following notation.  By 

\begin{equation}\label{item5sa}
(Y^{(1)}_n,Y^{(2)}_n,\dots,Y^{(k)}_n) \Rightarrow (Y^{(1)}_{\infty},Y^{(2)}_{\infty},\dots,Y^{(k)}_{\infty})
\end{equation}

\noindent we shall indicate the weak convergence of the 
{\it joint} law of the $k$-vector
$(Y^{(1)}_n,Y^{(2)}_n,\dots,Y^{(k)}_n)$ to that of
$(Y^{(1)}_{\infty},Y^{(2)}_{\infty},\dots,Y^{(k)}_{\infty})$,
as $n \rightarrow \infty$.
Since $LI_n$ is the length of the top row of the associated
Young diagrams, the length of the second row is simply $n - LI_n$. 
Denoting the length of the $i^{th}$ row by $R^{i}_n$, \eqref{item5sa},
together with an application of the Cram\'er-Wold Theorem,
recovers the result of Chistyakov and G\"otze \cite{ChG}
as part of the following easy corollary,
which is in fact equivalent to Theorem \ref{thm1}:\\

\begin{corollary} \label{cor1}
For the sequence in Theorem \ref{thm1},
if $a=b>0$, then

\begin{equation}\label{item5t}
\left(\frac{R^{1}_n-n/2}{\sqrt {n}}, \frac{R^{2}_n-n/2}{\sqrt{n}}\right)  
\Rightarrow  Y_{\infty} := (R^{1}_{\infty},R^{2}_{\infty}),
\end{equation}

\noindent where the law of $Y_{\infty}$ is 
supported on the $2^{nd}$ main diagonal of $\bbr^2$,
and with

$$R^{1}_{\infty} \stackrel{\cal{L}}{=}  \sqrt{\frac{1-a}{a}}\left( - \frac12 B(1)+ \max_{0\le t \le 1}B(t)\right).$$

\noindent If $a \ne b$ or $a=b=0$, then setting
$\pi_{min} = \min\{\pi_1,\pi_2\}$, we have

\begin{equation}\label{item5u}
\left(\frac{R^1_n  - n\pi_{max}}{\sqrt{n}} ,\frac{R^2_n  - n\pi_{min}}{\sqrt{n}}\right)
\Rightarrow   N((0,0),\tilde{\Sigma}),
\end{equation}

\noindent where $\tilde{\Sigma}$ is the covariance matrix

$$(\tilde{\sigma}^2/4) \begin{pmatrix} 1 &-1  \\ -1 & 1   \end{pmatrix},$$

\noindent where $\tilde{\sigma}^2 = 4ab(2-a-b)/(a+b)^3$, for $a \ne b$,
and $\tilde{\sigma}^2=0$, for $a=b=0$.

\end{corollary}

\begin{Rem}
The joint distributions in \eqref{item5t} and \eqref{item5u}
are of course degenerate, in  that the sum of the two
components is a.s.~identically zero in each case.  In \eqref{item5t},
the density of the first component of $R_{\infty}$ is
easy to find, and is given by
(e.g., see \cite{HLM})

\begin{equation}\label{item5ua}
f(y) = \frac{16}{\sqrt{2\pi}} \left(\frac{a}{1-a}\right)^{3/2} y^2 e^{-2ay^2/(1-a)}, \qquad y \ge 0.
\end{equation}

\noindent As in Chistyakov and G\"otze \cite{ChG}, 
\eqref{item5t} can then be stated as:
For any bounded, continuous function 
$g:\bbr^2 \rightarrow \bbr$,

\begin{align*}
&\lim_{n \rightarrow \infty} \left( g \left(\frac{R^1_n-n/2}{\sqrt {(1-a)n/a}}, 
\frac{R^2_n-n/2}{\sqrt{(1-a)n/a}}\right) \right)\nonumber\\
&\qquad = 2\sqrt{2\pi}\int_0^{\infty}  g(x,-x)\phi_{GUE,2}(x,-x) dx,
\end{align*}

\noindent where $\phi_{GUE,2}$ is the density of the
eigenvalues of the $2 \times 2$ GUE, and
is given by

$$\phi_{GUE,2}(x_1,x_2) = \frac{1}{\pi}(x_1-x_2)^2 e^{-(x_1^2 + x_2^2)}.$$

To see the GUE connection more explicitly, consider the  
$2 \times 2$ traceless GUE matrix

\begin{equation*}
M_0 =\begin{pmatrix}
X_1 & Y + iZ\\
Y - iZ  &X_2
\end{pmatrix},
\end{equation*}

\noindent where $X_1, X_2, Y$, and $Z$
are centered, normal random variables.
Since $\mbox{Corr }(X_1,X_2) = -1$, 
the largest eigenvalue of $M_0$ is

$$\lambda_{1,0} = \sqrt{X_1^2 + Y^2 + Z^2},$$

\noindent almost surely, 
so that $\lambda_{1,0}^2 \sim \chi_3^2$
if $\mbox{Var } X_1 =  \mbox{Var } Y = \mbox{Var } Z = 1$.
Hence, up to a scaling factor, the density of $\lambda_{1,0}$ is
given by \eqref{item5ua}.  
Next, let us perturb $M_0$ to

$$M = \alpha GI + \beta M_0,$$

\noindent where $\alpha$ and $\beta$ are constants, $G$
is a standard normal random variable independent of $M_0$,
and $I$ is the identity matrix.  
The covariance of the diagonal elements of $M$ is then computed to be
$\rho := \alpha^2 - \beta^2$.  Hence, to obtain a desired value
of $\rho$, we may take $\alpha = \sqrt{(1+\rho)/2}$ and
$\beta = \sqrt{(1-\rho)/2}.$  Clearly, the largest eigenvalue
of $M$ can then be expressed as

\begin{equation}\label{item5ub}
\lambda_1 = \sqrt{\frac{1+\rho}{2}}G + \sqrt{\frac{1-\rho}{2}}\lambda_{1,0}.
\end{equation}

\noindent At one extreme, $\rho = -1$, we recover
$\lambda_1 = \lambda_{1,0}$.  At the other extreme, $\rho = 1$,
we obtain $\lambda_1 = Z$.  Midway between these two extremes,
at $\rho=0$, we have a standard GUE matrix, so that

$$\lambda_1 =  \sqrt{\frac{1}{2}}\left(G + \lambda_{1,0}\right).$$
\end{Rem}

\end{document}